\documentclass[12pt]{amsart}
\usepackage{amsmath,amssymb,latexsym,amsthm}
\usepackage[small]{diagrams}
\diagramstyle[labelstyle=\scriptstyle,heads=littlevee]

\theoremstyle{plain}
\newtheorem{theorem}[subsection]{Theorem}
\newtheorem{proposition}[subsection]{Proposition}
\newtheorem{corollary}[subsection]{Corollary}
\newtheorem{lemma}[subsection]{Lemma}

\newcommand{\bksl}{{\smallsetminus}}

\newcommand{\ch}{{\operatorname{ch}}}
\newcommand{\conv}{{\operatorname{conv}}}
\newcommand{\Irr}{{\operatorname{Irr}}}

\newcommand{\bR}{{\mathbb{R}}}
\newcommand{\bZ}{{\mathbb{Z}}}

\newcommand{\Aut}{{\operatorname{Aut}}}

\newcommand{\vt}{{\operatorname{vert}}}
\newcommand{\edg}{{\operatorname{edg}}}

\newcommand{\polytope}{{\Pi}}
\newcommand{\edges}{{\operatorname{E}}}
\newcommand{\weights}{{\operatorname{wts}}}

\newcommand{\tK}{{\operatorname{K}}}
\newcommand{\tW}{{\operatorname{W}}}
\newcommand{\tX}{{\operatorname{X}}}
\newcommand{\tY}{{\operatorname{Y}}}

\newcommand{\SGAthreeIII}{{SGA\,{\mathsurround=0pt$3_{{}^{\text{III}}}$}}}

\newcommand{\fact}{{\smaller{\textsc{Fact}}}}

\begin{document}

\title[Determining a reductive group from its irreducible representations]
{Determining a connected split reductive group \\
from its irreducible representations}

\author{CheeWhye Chin}

\address{Department of Mathematics\\
         National University of Singapore\\
	 2 Science Drive 2\\
	 Singapore 117543\\
	 Singapore}


\email{cheewhye@math.nus.edu.sg}


\date{\today}




%
%

\begin{abstract}

We show that
a connected split reductive group $G$
over a field of characteristic~$0$
is uniquely determined
up to isomorphism
by specifying
a maximal torus $T$ of $G$,
the set of isomorphism classes
of irreducible representations of $G$,
and the character homomorphism
from the Grothendieck ring of $G$ 
to that of $T$.

\end{abstract}

\maketitle



\section{The main result}

\subsection{}
Let $G$ be
a connected split
reductive algebraic group
over a field $k$.
The finite dimensional
$k$-rational representations of $G$
will be simply called
representations of $G$.
The set $\Irr(G)$ of
isomorphism classes of
irreducible representations of $G$
forms a free $\bZ$-basis of
the Grothendieck ring $\tK(G)$ of $G$,
and one has
the canonical inclusion
\[
   \Irr(G)
   \rInto
   \tK(G)
   \qquad\text{of sets}.
\]
Let $T \subseteq G$ be
a maximal torus of $G$,
and let $\tX(T)$ denote
the free abelian group of
characters of $T$;
one then has
the canonical inclusion
\[
   \tX(T)
   \rInto
   \tK(T)
   \qquad\text{of multiplicative monoids}.
\]
Restricting representations
from $G$ to $T$
induces 
the ``character homomorphism''
\[
   \ch_G
   \ :\ 
   \tK(G)
   \rInto
   \tK(T)
   \qquad\text{of rings},
\]
which is injective
because $G$ is connected.

\subsection{}
Consider
the question of
determining
the group $G$
from the canonical maps
\[
   \Irr(G)
   \rInto
   \tK(G)
   \rInto^{\quad \ch_G \quad}
   \tK(T)
   \lInto
   \tX(T),
\]
which may be viewed as
an ``inverse problem''
in the representation theory
for connected split
reductive groups.
The main result
of this paper,
theorem~(\ref{thm:main}),
asserts that
the group $G$ is indeed
\textit{uniquely} determined
up to isomorphism
by the maps above.
It was worked out
for an application to
the question of
``independence of $\ell$ of monodromy groups''
but may be of
independent interest.

\subsection{Hypotheses for theorem~(\ref{thm:main})}
\label{hyp:thm:main}
Let $k$ be a field
of characteristic~$0$,
let $G$ and $G'$ be
connected split
reductive algebraic groups over $k$,
and let $T \subseteq G$
and $T' \subseteq G'$ be
maximal tori.
Suppose we are given
a bijection of sets
$ \phi
  :
  \Irr(G')
  \rTo^{\simeq}
  \Irr(G)
$
and an isomorphism of tori
$ f_T
  :
  T
  \rTo^{\simeq}
  T'
$
which are
compatible with each other
in the sense that
the diagram
\[
  \begin{diagram}
    \tX(T')
    &
    \rTo_{\quad \tX(f_T) \quad}^{\simeq}
    &
    \tX(T)
    \\
    \dInto
    &
    &
    \dInto
    \\
    \tK(T')
    &
    \rTo_{\quad \tK(f_T) \quad}^{\simeq}
    &
    \tK(T)
    \\
    \uInto^{\ch_{G'}}
    &
    &
    \uInto_{\ch_G}
    \\
    \tK(G')
    &
    \rTo_{\quad \tK(\phi) \quad}^{\simeq}
    &
    \tK(G)
    \\
    \uInto
    &
    &
    \uInto
    \\
    \Irr(G')
    &
    \rTo_{\quad \phi \quad}^{\simeq}
    &
    \Irr(G)
  \end{diagram}
  \qquad
  \text{commutes}.
\]
Here,
$\tX(f_T)$
(resp.~$\tK(f_T)$) is
the group isomorphism
(resp.~ring isomorphism)
induced by $f_T$,
and $\tK(\phi)$ is
the isomorphism of additive groups
induced by $\phi$.
Note that
since the ring homomorphisms
$\ch_G$ and $\ch_{G'}$
are injective,
the assumption that
the above diagram commutes
implies that
$\tK(\phi)$ is necessarily
a ring isomorphism.

\begin{theorem}
\label{thm:main}
Under the hypotheses
of~(\ref{hyp:thm:main}),
the following conclusions hold.
\begin{itemize}
\item[(a)]
There exists
an isomorphism of algebraic groups
$ f
  :
  G
  \rTo^{\simeq}
  G'
$
compatible with
$f_T$ and $\phi$
in the sense that
\begin{itemize}
\item[(i)]
$f$ extends
the isomorphism $f_T$
of maximal tori,
i.e.~the diagram
\[
  \begin{diagram}
    T
    &
    \rTo^{\quad f_T \quad}_{\eqsim}
    &
    T'
    \\
    \dInto
    &
    &
    \dInto
    \\
    G
    &
    \rTo_{\quad f \quad}^{\simeq}
    &
    G'
  \end{diagram}
  \qquad\text{commutes};
  \quad\text{and}
\]
\item[(ii)]
the bijection of sets
\[
  \begin{aligned}
    f^{*}
    &
    \ :\ 
    \Irr(G')
    \rTo^{\simeq}
    \Irr(G)
    \qquad
    \text{given by
          $\rho' \mapsto \rho'\cdot f$}
    \\
    \text{is equal to}
    \qquad
    \phi
    &
    \ :\ 
    \Irr(G')
    \rTo^{\simeq}
    \Irr(G).
  \end{aligned}
\]
\end{itemize}
\item[(b)]
If
$ \widetilde{f}
  :
  G
  \rTo^{\simeq}
  G'
$
is another isomorphism of algebraic groups
compatible with
$f_T$ and $\phi$,
then there exists
a $k$-rational point
$t \in T(k)$ of $T$
such that
\[
   \widetilde{f}
   \ =\ 
   f \cdot \text{(conjugation by $t$)}.
\]
\end{itemize}
\end{theorem}

\subsection{}
Since
the structure theory
of a connected split
reductive group
is very closely related to
its representation theory,
theorem~(\ref{thm:main})
should be
``almost obvious''
to the experts.
Indeed,
one just has to show that
the hypotheses
of~(\ref{hyp:thm:main})
induce
an isomorphism between
the root data of
$G$ and $G'$
and apply
the ``isomorphism theorem''
to conclude that
$G$ and $G'$
are isomorphic.
Let us
make this argument
more precise.
Recall that
$\tX(T)$
denotes
the group of characters of $T$;
let $\tY(T)$ be
the group of cocharacters of $T$,
and let $\Phi(G,T)$ (resp.~$\Phi^\vee(G,T)$)
denote
the set of roots (resp.~coroots)
of $G$ with respect to $T$.
The given isomorphism
$ f_T
  :
  T
  \rTo^{\simeq}
  T'
$
of maximal tori
already induces
an adjoint pair of
isomorphisms
\[
   \tX(f_T)
   \ :\ 
   \tX(T')
   \rTo^{\simeq}
   \tX(T)
   \qquad
   \text{and}
   \qquad
   \tY(f_T)
   \ :\ 
   \tY(T)
   \rTo^{\simeq}
   \tY(T')
\]
on the groups of
characters
and
cocharacters
of the tori,
so it is a matter
of showing that
$\tX(f_T)$ maps
$\Phi(G',T')$
onto
$\Phi(G,T)$
and that
$\tY(f_T)$ maps
$\Phi^\vee(G,T)$
onto
$\Phi^\vee(G',T')$.

What is not clear at all
from the hypotheses
of~(\ref{hyp:thm:main})
is that
$\tX(f_T)$ should send
$\Phi(G',T')$
\textit{into} $\Phi(G,T)$.
Indeed,
suppose $\rho'$ is
the adjoint representation
of $G'$,
and to fix ideas,
let us furthermore
assume that
$\rho'$ is
irreducible.
The non-zero weights
of $\rho'$
are then
the roots $\Phi(G',T')$,
and it is easy to deduce
from the hypotheses
of~(\ref{hyp:thm:main})
that
$\tX(f_T)$ maps these
to the non-zero weights 
of the irreducible representation
$\rho := \phi(\rho')$
of $G$.
If we somehow know that
$\rho$ is also
the adjoint representation
of $G$,
then it would follow that
$\tX(f_T)$
sends $\Phi(G',T')$
onto $\Phi(G,T)$;
however,
there is no
\textit{a priori} reason
to assume that
$\rho$ has anything to do
with the adjoint representation of $G$.
It is therefore necessary to
justify the validity of
the following:

\begin{proposition}
\label{bijection:roots}
Under the hypotheses
of~(\ref{hyp:thm:main}),
the isomorphism
\[
  \tX(f_T)
  \ :\ 
  \tX(T')
  \rTo^{\simeq}
  \tX(T)
  \qquad
  \text{maps
        $\Phi(G',T')$
        onto
	$\Phi(G,T)$}.
\]
\end{proposition}

We will establish this proposition
in \S\ref{sect:weight polytopes and roots};
for now,
let us proceed
to deduce theorem~(\ref{thm:main})
from it.
One knows that
a subset
$\Delta'\subseteq\Phi(G',T')$
is a system of simple roots
if and only if
every element of $\Phi(G',T')$
is a unique $\bZ$-linear combination
of elements in $\Delta'$
with coefficients of the same sign.
Hence
(\ref{bijection:roots})
and the linearity of $\tX(f_T)$
imply:

\begin{corollary}
\label{bijection:simple roots}
Under the hypotheses
of~(\ref{hyp:thm:main}),
the isomorphism
\[
  \tX(f_T)
  \ :\ 
  \tX(T')
  \rTo^{\simeq}
  \tX(T)
\]
maps
a system of simple roots of $(G',T')$
to
a system of simple roots of $(G,T)$.
\end{corollary}

Next,
consider the Weyl group
$W(G',T') \subseteq \Aut_\bZ(\tX(T'))$.
Since $\tK(T')$ is
the group ring
over $\tX(T')$,
the group
$\Aut_\bZ(\tX(T'))$
and hence $W(G',T')$
acts faithfully on
the ring $\tK(T')$
by ring automorphisms.
One knows that
the subring of
$W(G',T')$-invariants in $\tK(T')$
is precisely $\tK(G')$.
It then follows
from Galois theory
that
$W(G',T')$
is identified with
the group of automorphisms of
the algebra $\tK(T')$
over $\tK(G')$.
Hence the commutativity of
the diagram in~(\ref{hyp:thm:main})
implies that:

\begin{proposition}
\label{bijection:Weyl group}
Under the hypotheses
of~(\ref{hyp:thm:main}),
the isomorphism
\[
   [\tX(f_T)]
   \ :\ 
   \Aut_\bZ(\tX(T'))
   \rTo^{\simeq}
   \Aut_\bZ(\tX(T))
   \qquad
   \text{given by}
   \qquad
   \sigma
   \rMapsto
   \tX(f_T) \cdot \sigma \cdot \tX(f_T)^{-1}
\]
maps
$\tW(G',T')$
onto
$\tW(G,T)$.
\end{proposition}

If $\alpha \in \Phi(G,T)$ is
a root,
write
$\alpha^\vee \in \Phi^\vee(G,T)$
for the coroot
determined by $\alpha$,
and write
$s_\alpha \in \Aut_\bZ(\tX(T))$ for
the corresponding
reflection
$ x
  \mapsto
  x
  -
  \langle x,\alpha^\vee \rangle\,\alpha
$.
One knows that
for any system
$\Delta\subseteq\Phi(G,T)$
of simple roots,
$s_\alpha$ is
the unique element of $\tW(G,T)$
which maps $\alpha$ to $-\alpha$
and stabilizes the subset
$P(\Delta)\bksl\{\alpha\}$,
where $P(\Delta)$ denotes
the subset of positive roots
determined by $\Delta$.
This fact,
together with~(\ref{bijection:simple roots})
and~(\ref{bijection:Weyl group}),
shows that:

\begin{corollary}
\label{bijection:reflections}
Under the hypotheses
of~(\ref{hyp:thm:main}),
if
$\alpha' \in \Phi(G',T')$,
and
$ \alpha
  :=
  \tX(f_T)(\alpha')
$,
then
\[
   [\tX(f_T)](s_{\alpha'})
   \ =\ 
   s_\alpha
   \qquad
   \text{in $\Aut_\bZ(\tX(T))$}.
\]
\end{corollary}

It is easy to see that
the coroot $\alpha^\vee\in\Phi^\vee(G,T)$
is the unique element $y\in\tY(T)$
with the property:
\[
  \text{for any $x\in\tX(T)$},
  \qquad
  \text{one has}
  \quad
  x - s_\alpha(x)
  \ =\ 
  \langle x,y \rangle \,\alpha
  \quad
  \text{in $\tX(T)$}.
\]
Thus~(\ref{bijection:reflections})
and the fact that
$ \tY(f_T) $
is the adjoint of
$ \tX(f_T) $
shows that
$\tY(f_T)(\alpha^\vee)$
is the unique element
$y' \in \tY(T')$
with the property:
\[
  \text{for any $x'\in\tX(T')$},
  \qquad
  \text{one has}
  \quad
  x' - s_{\alpha'}(x')
  \ =\ 
  \langle x',y' \rangle \,\alpha'
  \quad
  \text{in $\tX(T')$},
\]
and hence:

\begin{corollary}
\label{bijection:coroots}
Under the hypotheses
of~(\ref{hyp:thm:main}),
the isomorphism
\[
  \tY(f_T)
  \ :\ 
  \tY(T)
  \rTo^{\simeq}
  \tY(T')
  \qquad
  \text{maps
        $\Phi^\vee(G,T)$
        onto
	$\Phi^\vee(G',T')$}.
\]
More precisely,
if $\alpha \in \Phi(G,T)$,
and
$ \alpha'
  :=
  \tX(f_T)^{-1}(\alpha)
$,
then
\[
  \tY(f_T)(\alpha^\vee)
  \ =\ 
  (\alpha')^\vee
  \qquad
  \text{in $\Phi^\vee(G',T')$}.
\]
\end{corollary}

\subsection{Proof of theorem~(\ref{thm:main})
            assuming~(\ref{bijection:roots})}
\label{proof:thm:main}
%
%
%
Applying~(\ref{bijection:roots})
and~(\ref{bijection:coroots}),
we see that
$\tX(f_T)$ and its adjoint $\tY(f_T)$ induce
an isomorphism of root data
\[
   F
   \ :\ 
   \Bigl(\ 
     \tX(T'),\Phi(G',T'),
     \tY(T'),\Phi^\vee(G',T')
   \ \Bigr)
   \rTo^{\simeq}
   \Bigl(\ 
     \tX(T),\Phi(G,T),
     \tY(T),\Phi^\vee(G,T)
   \ \Bigr)
   .
\]
The isomorphism theorem
for connected split reductive groups
(cf.~\cite{SGAthreeIII} exp.~XXIII, Th.~4.1)
then shows that
the given isomorphism
$ f_T
  :
  T
  \rTo^{\simeq}
  T'
$
of tori
extends to
an isomorphism
$ f
  :
  G
  \rTo^{\simeq}
  G'
$
of groups
compatible with
the isomorphism~$F$
of root data,
and that
such an extension is unique
up to precomposition with
conjugation by
a $k$-rational point of $T$.
To complete
part~(a) of theorem~(\ref{thm:main}),
we have to verify that
the two maps
on irreducible representations
\[
  f^{*}
  \ :\ 
  \Irr(G')
  \rTo^{\simeq}
  \Irr(G)
  \qquad
  \text{and}
  \qquad
  \phi
  \ :\ 
  \Irr(G')
  \rTo^{\simeq}
  \Irr(G)
  \qquad
  \text{are equal}.
\]
Let $\rho'\in\Irr(G')$ be
an irreducible representation of $G'$,
and let
$\weights(\rho') \subseteq \tX(T')$ denote
its set of weights.
To show that
$\phi(\rho')$
and
$f^{*}(\rho')$
are equal in $\Irr(G)$,
it suffices for us to show that
their sets of weights
$\weights(\phi(\rho'))$
and $\weights(f^{*}(\rho'))$
are both equal to
$\tX(f_T)(\weights(\rho'))$
in $\tX(T)$.
The fact that
$ \weights(\phi(\rho'))
  =
  \tX(f_T)(\weights(\rho'))
$
follows from
the commutativity of
the diagram in~(\ref{hyp:thm:main}).
On the other hand,
the irreducible representation
$f^{*}(\rho') \in \Irr(G)$
of $G$
is obtained by
precomposition:
$ f^{*}(\rho')
  =
  \rho'\cdot f
$,
so
$ f^{*}(\rho')|_T
  =
  \rho'|_{T'}\cdot f_T
$
because $f$ restricts to $f_T$;
this implies that
$x'\in\tX(T')$ is a weight of $\rho'$
if and only if
$ \tX(f_T)(x')
  =
  x'\cdot f_T
  \in\tX(T)$
is a weight of $f^{*}(\rho')$,
and hence
$ \weights(f^{*}(\rho'))
  =
  \tX(f_T)(\weights(\rho')) $.
We have proved
part~(a) of theorem~(\ref{thm:main}).

To show
part~(b) of theorem~(\ref{thm:main}),
let
$ \widetilde{f}
  :
  G
  \rTo^{\simeq}
  G'
$
be any isomorphism of groups
having the same properties~(i) and~(ii)
as $f$ in part~(a).
Then
the hypotheses of~(\ref{hyp:thm:main})
are satisfied
with $\phi = \Irr(\widetilde{f})$
and $f_T = \widetilde{f}|_T$,
so~(\ref{bijection:roots})
and~(\ref{bijection:coroots})
now show that
$\tX(\widetilde{f}|_T)$
and its adjoint
$\tY(\widetilde{f}|_T)$
induce an
isomorphism
from the root datum
of $(G',T')$
to that
of $(G,T)$.
But this isomorphism of root data
is the same as~$F$
because $\tX(\widetilde{f}|_T) = \tX(f_T)$
and $\tY(\widetilde{f}|_T) = \tY(f_T)$.
As we have observed before,
it is then a consequence of
the isomorphism theorem for
connected split reductive groups that
\[
   \widetilde{f}
   \ =\ 
   f \cdot (\text{conjugation by $t$})
   \qquad
   \text{for some $k$-rational point
         $t\in T(k)$ of $T$}.
\]
\qed

\subsection{Remark}
I do not know
if theorem~(\ref{thm:main}) holds
when $k$ has characteristic~$>0$.
In the key steps
(\ref{bijection:roots})
--
(\ref{bijection:coroots})
leading to
the above proof~(\ref{proof:thm:main}),
the hypothesis that
$k$ is of characteristic~$0$
is needed only for
the proof of~(\ref{bijection:roots})
in \S\ref{sect:weight polytopes and roots}
below;
if~(\ref{bijection:roots}) holds
in positive characteristic,
the proofs of the other propositions
as well as
the argument in~(\ref{proof:thm:main})
go through without change.
%


\section{Weight polytopes and roots}
\label{sect:weight polytopes and roots}

Our goal now
is to establish
proposition~(\ref{bijection:roots}).
The main idea
is to describe
the roots of
a connected split
reductive group
in terms of
the geometry of
the weight polytopes
associated to
its irreducible
representations.
Let us first
fix some notation.

\subsection{}
\label{subsect:notation}
Let $G$ be
a connected split reductive group
over a field $k$,
let $T \subseteq G$ be
a maximal torus,
and let $\Phi(G,T)$ be
the set of roots of $G$
with respect to $T$.
If $\Delta\subseteq\Phi(G,T)$ is
a system of simple roots,
let
$ P(\Delta)
  :=
  \bR_{\geqslant 0}\,\Delta
  \cap
  \Phi(G,T) $
denote
the system of positive roots
determined by $\Delta$,
and let
\[
    C(\Delta)
    \ :=\ 
    \{\ 
        x \in \tX(T)_\bR
	\ :\ 
	\langle x,\alpha^\vee \rangle \geqslant 0
	\quad
	\text{for every}
	\quad
	\alpha \in \Delta
    \ \}
\]
denote
the dominant
closed Weyl chamber
determined by $\Delta$
(equivalently,
 by $P(\Delta)$);
it is a fundamental domain
for the action of
the Weyl group
$\tW(G,T)$
on the vector space
$\tX(T)_\bR$.
For any root $\alpha \in \Phi(G,T)$,
let
$ s_\alpha
  :=
  (
  x
  \mapsto
  x
  -
  \langle x,\alpha^\vee \rangle\,\alpha
  )
$
denote
the corresponding reflection
in the vector space $\tX(T)_\bR$.
If $x_0, x_1 \in \tX(T)_\bR$
are two points,
let $[x_0,x_1]$ denote
the convex hull
of $\{x_0,x_1\}$.
If $K \subseteq \tX(T)_\bR$
is a subset,
let us say that
$x \in \tX(T)_\bR$
is an \textsl{indivisible element of $K$}
if and only if
$x \in K$
and
for every integer $n > 1$,
$\tfrac{1}{n}\,x \notin K$.

If $\rho\in\Irr(G)$ is
an irreducible representation
of $G$,
its \textsl{weight polytope}
$\polytope(\rho)$
is defined as
the convex hull
of the set of weights
$\weights(\rho)$ of $\rho$
in the $\bR$-vector space
$\tX(T)_\bR$.
The action of $\tW(G,T)$
on $\tX(T)_\bR$
stabilizes
$\polytope(\rho)$
since it stabilizes
$\weights(\rho)$.
It is a well-known fact that 
the set of vertices
of the weight polytope
$\polytope(\rho)$
is precisely
the Weyl group orbit of
the highest weight
(with respect to
any system of simple roots)
of $\rho$.
It seems to be
a lesser-known fact that
the edges
of $\polytope(\rho)$
can also be described
in terms of
the root datum
of the pair $(G,T)$
and
the highest weight of $\rho$:

\begin{proposition}
\label{prop:wp edges}
Let $\rho \in \Irr(G)$ be
an irreducible representation of $G$,
let $\Delta\subseteq\Phi(G,T)$ be
a system of simple roots,
and let $x_0 \in \tX(T)$ be
the highest weight of $\rho$
with respect to $\Delta$.
Let $\edges(\rho,{x_0})$ denote
the set of edges of $\polytope(\rho)$
emanating from
the vertex $x_0$.
Then
there exists
a subset
$ \Delta(x_0) \subseteq P(\Delta)$
of positive roots
and a bijection
\[
   \Delta(x_0)
   \rTo^{\simeq}
   \edges(\rho,{x_0}),
   \qquad
   \alpha
   \mapsto
   [x_0,s_\alpha(x_0)].
\]
Furthermore,
if $x_0$ is
strongly dominant
with respect to $\Delta$,
then
$\Delta(x_0) = \Delta$.
\end{proposition}

In other words,
every edge of $\polytope(\rho)$
emanating from $x_0$
is obtained
by reflecting $x_0$
in the wall
of a uniquely determined
positive root $\alpha$,
and when $x_0$ is
strongly dominant,
the possible $\alpha$'s
are precisely
the simple roots.
This result is
intuitively clear
from the geometry of
the weight polytopes;
anyhow,
we will prove
a more precise
theorem~(\ref{thm:wp edges})
in \S\ref{sect:wp edges}.
For showing
proposition~(\ref{bijection:roots}),
we will need to 
recover the root $\alpha$
from the edge
$[x_0,s_\alpha(x_0)]$
it determines,
and this is given by
the following:

\begin{proposition}
\label{prop:root from edge}
Suppose that
the base field $k$
is of characteristic~$0$.
In the situation of
proposition~(\ref{prop:wp edges}),
let $e \in E(\rho,x_0)$
be an edge of $\polytope(\rho)$
emanating from
the vertex $x_0$;
then
the positive root
$\alpha \in \Delta(x_0)$
for which
$e = [x_0,s_\alpha(x_0)]$
is characterized as
the unique indivisible element of
$x_0 - \bigl( \weights(\rho) \cap e \bigr)$.
\end{proposition}

\begin{proof}
Let
$ c
  :=
  \langle x_0,\alpha^\vee \rangle
  \in \bZ$;
this is an integer $>0$
since $x_0\in\tX(T)$ is
dominant with respect to $\Delta$
and
$s_\alpha(x_0) = x_0 - c\,\alpha$
is $\ne x_0$.
It suffices for us
to show that
\[
   \weights(\rho) \cap e
   \ =\ 
   \Bigl\{\ 
      x_0 - t\,\alpha
      \in \tX(T)
      \ :\ 
      t \in \{ 0,1,\ldots,c \}
   \ \Bigr\}
   \qquad
   \text{in $\tX(T)_\bR$}.
\]
Suppose $x \in \weights(\rho) \cap e$.
Then
$x_0 - x$ lies in
$\bZ\,\Delta$
because
$x_0$ is the highest weight of $\rho$,
and it lies in
$[0,c\,\alpha]$
because
$x_0 - e = [0,c\,\alpha]$.
So to show the inclusion $\subseteq$,
we just have to show that
$ \bZ\,\Delta \cap [0,c\,\alpha]
  =
  \{ 0,\alpha,\ldots,c\,\alpha \}$.
Pick a Weyl group element
$w\in\tW(G,T)$
such that $w(\Delta)$ contains $\alpha$.
Then
$ \bZ\,\Delta \cap [0,c\,\alpha]
  =
  \bZ\,w(\Delta) \cap [0,c\,\alpha]$,
and this is equal to
$\{ 0,\alpha,\ldots,c\,\alpha \}$
because 
$w(\Delta)$ is
a linearly independent subset
of the $\bR$-vector space $\tX(T)_\bR$
and it contains $\alpha$.

For the reverse inclusion $\supseteq$,
we have to show that 
$x_0 - t\,\alpha$
for $t = 0,1,\ldots,c$
actually occur
as weights of $\rho$.
By a standard argument
(see the proof of
 \cite{Serre-GroupesGrothendieck}~Lemme~5
 for instance),
we are reduced to the case
when the base field $k$
is algebraically closed
and the group $G$
is a connected semisimple group.
The assertion we want
is then
a consequence of
the corresponding assertion of
the representation of the Lie~algebra of $G$
(cf.~\cite{Samelson-NotesLieAlg}
 \S3.2 Th.~B(c) for instance);
here, we have used
the assumption that
$k$ is of characteristic~$0$.
\end{proof}

\subsection{Proof of proposition~(\ref{bijection:roots})}
We assume
the hypotheses of~(\ref{hyp:thm:main}).
By symmetry,
it suffices
for us to show that
$\tX(f_T)$ maps
$\Phi(G',T')$
\textit{into}
$\Phi(G,T)$.
Let
$\alpha' \in \Phi(G',T')$ be given
and let
$\alpha := \tX(f_T)(\alpha') \in\tX(T)$;
we want to show that
$\alpha$ is in fact a root of $(G,T)$.

Choose a system
$\Delta' \subseteq\Phi(G',T')$
of simple roots of $(G',T')$
containing $\alpha'$,
and choose
a character
$x'_0 \in \tX(T') \cap C(\Delta')$
of $T'$
\textit{strongly dominant}
with respect to $\Delta'$.
Let $\rho' \in \Irr(G')$ be
the irreducible representation of $G'$
have $x'_0$ as
highest weight with respect to $\Delta'$.
By~(\ref{prop:wp edges}),
$ e' := [x'_0,s_{\alpha'}(x'_0)] $
is an edge of
the weight polytope
$\polytope(\rho')$ of $\rho'$
emanating from
the vertex $x'_0$,
and
by~(\ref{prop:root from edge}),
\[
  \alpha'
  \quad
  \text{is the unique indivisible element of}
  \quad
  x'_0 - \big( \weights(\rho') \cap e' \big).
\]
Let
$\rho := \phi(\rho') \in\Irr(G)$.
The commutativity of
the diagram in~(\ref{hyp:thm:main})
shows that
$\tX(f_T)$ maps
$\weights(\rho') \subseteq \tX(T')$
onto
$\weights(\rho) \subseteq \tX(T)$.
The linearity of $\tX(f_T)$
then ensures that
the polytope
$\polytope(\rho') \subseteq \tX(T')_\bR$
is mapped onto
the polytope
$\polytope(\rho) \subseteq \tX(T)_\bR$.
It follows that
$ e := \tX(f_T)(e') $
is an edge of $\polytope(\rho)$
emanating from
the vertex
$x_0 := \tX(f_T)(x'_0)$,
and that
\[
   \alpha
   \ :=\  
   \tX(f_T)(\alpha')
   \quad
   \text{is the unique indivisible element of}
   \quad
   x_0 - \big( \weights(\rho) \cap e \big).
\]
Now choose a system
$\Delta \subseteq \Phi(G,T)$
of simple roots of $(G,T)$
so that
$x_0$ is dominant
with respect to $\Delta$
\footnote{\ 
Note that
there is no \textit{a priori} reason
to assume that
$x_0$ is \textit{strongly dominant}
with respect to
any system of
simple roots of $(G,T)$,
even though
$x_0$ is the image of $x_0'$
which is strongly dominant
with respect to $\Delta'$.}.
It is then clear that
$x_0$ is
the highest weight of $\rho$
with respect to $\Delta$,
since $x_0$ lies in $C(\Delta)$
and is a vertex
of $\polytope(\rho)$.
We can now
apply~(\ref{prop:wp edges}) again
to see that
the edge $e$
is of the form
$[x_0, s_\beta(x_0)] $
for some root $\beta\in\Phi(G,T)$
which is positive
with respect to $\Delta$,
and by~(\ref{prop:root from edge}),
\[
   \beta
   \quad
   \text{is also the unique indivisible element of}
   \quad
   x_0 - \big( \weights(\rho) \cap e \big).
\]
Therefore
$\alpha = \beta$
lies in $\Phi(G,T)$.
\qed


\section{Edges of weight polytopes}
\label{sect:wp edges}

In this section,
we prove
theorem~(\ref{thm:wp edges})
which is a more precise version of
proposition~(\ref{prop:wp edges}).
To state the result,
let us fix
the necessary notation.

\subsection{}
\label{wp edges notation}
We keep the notation
introduced in~(\ref{subsect:notation}).
In addition,
let $\rho \in \Irr(G)$ be
an irreducible representation of $G$,
let $\Delta\subseteq\Phi(G,T)$ be
a system of simple roots,
and let $x_0 \in \tX(T)$ be
the highest weight of $\rho$
with respect to $\Delta$.
Define
\[
  \Phi_0
  \ :=\ 
  \{\ 
      \beta \in \Phi(G,T)
      \ :\ 
      \langle x_0,\beta^\vee \rangle = 0
  \ \},
  \qquad
  P_0
  \ :=\ 
  \Phi_0 \cap P(\Delta),
  \qquad
  \Delta_0
  \ :=\ 
  \Phi_0 \cap \Delta,
\]
let
$\tW_0$ be
the subgroup of $\tW(G,T)$
generated by
$ \{\ 
      s_\beta \in \tW(G,T)
      \ :\ 
      \beta \in \Delta_0
  \ \}$,
and let
$ \Delta(x_0)
  \ :=\ 
  \tW_0\,(\Delta\bksl\Delta_0)
  \ \subseteq
  \Phi(G,T) 
$
be the union of
the transforms of
$\Delta\bksl\Delta_0$
under $\tW_0$.
Finally,
let $\edges(\rho,{x_0})$ denote
the set of edges of $\polytope(\rho)$
emanating from
the vertex $x_0$.

\begin{theorem}
\label{thm:wp edges}
With the notation
of~(\ref{wp edges notation}),
one has
$ \Delta\bksl\Delta_0
  \ \subseteq\ 
  \Delta(x_0)
  \ \subseteq\  
  P(\Delta)\bksl P_0$,
and
the map
\[
   \Delta(x_0)
   \rTo
   \edges(\rho,{x_0}),
   \qquad
   \alpha
   \mapsto
   [x_0,s_\alpha(x_0)],
   \qquad
   \text{is a well-defined bijection}.
   \tag{*}
   \label{thm:wp edges:bijection}
\]
\end{theorem}

To illustrate:
if $x_0\in\tX(T)$ is
\textit{strongly dominant} with respect to $\Delta$,
then
the subset
$\Delta_0 \subseteq \Delta$ above
is the empty set,
$\tW_0$ is
the trivial subgroup of $\tW(G,T)$,
and $\Delta(x_0)$ is equal to $\Delta$.
On the other hand,
if $\rho$ is
a one-dimensional representation of $G$,
then $\Delta_0 = \Delta$,
$\tW_0 = \tW(G,T)$,
and $\Delta(x_0)$ is the empty set
---
i.e.~the polytope $\polytope(\rho)$
reduces to a singleton $\{x_0\}$
and hence has no edges.

We first establish
a series of lemmas
necessary for the proof.
The hypotheses
of the theorem
are assumed
throughout
the rest of this section.

\begin{lemma}
\label{wp edges:lemma1}
One has
$ \Delta\bksl\Delta_0
  \ \subseteq\ 
  \Delta(x_0)
  \ \subseteq\  
  P(\Delta)\bksl P_0$.
\end{lemma}

\begin{proof}
The first inclusion
is clear.
From the definitions
it is also clear that
$\Delta\bksl\Delta_0$
is contained in
$P(\Delta)\bksl P_0$.
Since
$\Delta(x_0) = \tW_0\,(\Delta\bksl\Delta_0)$,
the second inclusion would follow
if we show that
$P(\Delta)\bksl P_0$
is stable under $\tW_0$.
Now $\tW_0$ is generated by 
the $s_\beta$ with $\beta\in \Delta_0\subseteq P_0$,
so it suffices for us to show that
$P(\Delta)\bksl P_0$
is stable under $s_\beta$
for any $\beta\in P_0$.
But if $\beta\in P_0$,
we can write
$P(\Delta)\bksl P_0
 =
 (P(\Delta) \bksl \{\beta\}) \bksl \Phi_0$,
and we just have to note that
both $P(\Delta)\bksl\{\beta\}$
and $\Phi_0$
are individually stable under $s_\beta$.
\end{proof}

\begin{lemma}
\label{wp edges:lemma2:step1}
One has
\[
  \bigcup_{w\in \tW_0}
     \{\ 
         x \in \tX(T)_\bR
	 \ :\ 
	 \langle x,\beta^\vee \rangle \geqslant 0
	 \quad
	 \text{for every}
	 \quad
	 \beta\in w(P_0)
     \ \}
  \ =\ 
  \tX(T)_\bR.
\]
\end{lemma}

\begin{proof}
Since
$\bZ_{>0}\Phi_0 \cap \Phi(G,T)
 \subseteq
 \Phi_0$
and
$\Phi_0 = -\Phi_0$,
it follows that
$(\tX(T),\Phi_0,\tY(T),\Phi_0^\vee)$
is itself a root datum
(cf.~\cite{SGAthreeIII} exp.~XXI, Prop.~3.4.1),
with
$\Delta_0 \subseteq \Delta$ as
a system of simple roots,
with $P_0 \subseteq P(\Delta)$ as
the corresponding
system of positive roots,
and $\tW_0 \subseteq\tW(G,T)$ as
its Weyl group.
If
\[
   D
   \ :=\ 
   \{\ 
        x \in \tX(T)_\bR
	\ :\ 
	\langle x,\beta^\vee \rangle \geqslant 0
	\quad
	\text{for every}
	\quad
	\beta\in P_0
   \ \}
\]
is the closed Weyl chamber
in $\tX(T)_\bR$
determined by $P_0$,
its image
under the action of $w^{-1}$
for any $w\in \tW_0$
admits the description
\[
   w^{-1}(D)
   \ =\ 
   \{\ 
       x \in \tX(T)_\bR
       \ :\ 
       \langle x,\beta^\vee \rangle \geqslant 0
       \quad
       \text{for every}
       \quad
       \beta\in w(P_0)
   \ \}.
\]
Since $D$ is
a fundamental domain
for the action of $\tW_0$
on $\tX(T)_\bR$,
the lemma follows.
\end{proof}

\begin{lemma}
\label{wp edges:lemma2:step2}
Let
$ C_0
  :=
  \tW_0\,C(\Delta)
  =
  \bigcup_{w\in \tW_0}
     w(C(\Delta))
$
be
the union of 
the transforms of $C(\Delta)$
under $\tW_0$.
Then
\[
   C_0
   \ =\ 
   \{\ 
       x \in \tX(T)_\bR
       \ :\ 
       \langle x,\alpha^\vee \rangle \geqslant 0
       \quad
       \text{for every}
       \quad
       \alpha\in \Delta(x_0)
   \ \}.
\]
\end{lemma}

\begin{proof}
The inclusion relations
in~(\ref{wp edges:lemma1})
yields
$ \Delta
  \subseteq
  \Delta(x_0)\cup P_0
  \subseteq
  P(\Delta)
$.
In view of
the definition~(\ref{subsect:notation})
of $C(\Delta)$,
we may write
\[
   C(\Delta)
   \ =\ 
   \{\ 
       x \in \tX(T)_\bR
       \ :\ 
       \langle x,\alpha^\vee \rangle \geqslant 0
       \quad
       \text{for every}
       \quad
       \alpha\in \Delta(x_0) \cup P_0
   \ \}.
\]
Therefore,
\[
   C_0 
   \ =\ 
   \tW_0\,C(\Delta)
   \ =\ 
   \bigcup_{w\in \tW_0}
      \{\ 
          x \in \tX(T)_\bR
	  \ :\ 
	  \langle x,\alpha^\vee \rangle \geqslant 0
	  \quad
	  \text{for every}
	  \quad
	  \alpha\in w(\Delta(x_0) \cup P_0)
      \ \}.
\]
But for every $w\in \tW_0$,
we have
\[
   w(\Delta(x_0) \cup P_0)
   \ =\ 
   w(\Delta(x_0)) \ \cup\  w(P_0)
   \ =\ 
   \Delta(x_0)\ \cup\  w(P_0),
\]
so we can write
$C_0 = C_1 \cap C_2$
where
\[
   \begin{aligned}
   C_1
   \ :=\ 
   &
   \{\ 
       x \in \tX(T)_\bR
       \ :\ 
       \langle x,\alpha^\vee \rangle \geqslant 0
       \quad
       \text{for every}
       \quad
       \alpha\in \Delta(x_0)
   \ \},
   \\
   C_2
   \ :=\ 
   &
   \bigcup_{w\in \tW_0}
      \{\ 
          x \in \tX(T)_\bR
	  \ :\ 
	  \langle x,\beta^\vee \rangle \geqslant 0
	  \quad
	  \text{for every}
	  \quad
	  \beta\in w(P_0)
      \ \}.
   \end{aligned}
\]
But
$C_2 = \tX(T)_\bR$
by~(\ref{wp edges:lemma2:step1}),
so $C_0 = C_1$,
and this proves the lemma.
\end{proof}

\begin{lemma}
\label{wp edges:lemma2}
Every point in
the Weyl group orbit
$\tW(G,T)\, x_0$ of $x_0$
is of the form
\[
   x_0
   -
   \sum_{\alpha\in \Delta(x_0)}
      c_\alpha\,\alpha
   \qquad
   \text{for some coefficients
         $c_\alpha\in\bZ_{\geqslant 0}$,
	 $\alpha\in \Delta(x_0)$}.
\]
\end{lemma}

\begin{proof}
We deduce this from~(\ref{wp edges:lemma2:step2})
following the usual argument
for the case when $\Delta(x_0) = \Delta$
(cf.~\cite{Samelson-NotesLieAlg} \S2.11 Lemma~J).
Define a partial order relation $\succcurlyeq$
in $\tX(T)$
as follows:
for any $x,x_1 \in\tX(T)$,
the relation $x_1 \succcurlyeq x$ holds
if and only if
\[
   x_1
   \ =\ 
   x
   +
   \sum_{\alpha\in \Delta(x_0)}
      c_\alpha\,\alpha
   \qquad
   \text{for some coefficients
         $c_\alpha\in\bZ_{\geqslant 0}$,
	 $\alpha\in \Delta(x_0)$}.
\]
Given a point
$x\in \tW(G,T)\, x_0$
in the Weyl group orbit of $x_0$,
we need to show that
$x_0 \succcurlyeq x$.
Choose an element
$x_1 \in \tW(G,T)\, x_0$
which is $\succcurlyeq x$
and is maximal
with respect to $\succcurlyeq$.
If $x_1$ does not lie in
the set $C_0$ of~(\ref{wp edges:lemma2:step2}),
then there would be
some $\alpha \in \Delta(x_0)$
such that
$\langle x_1,\alpha^\vee \rangle$
is $< 0$;
but this means that
$ s_\alpha(x_1)
  =
  x_1
  -
  \langle x_1,\alpha^\vee \rangle \,\alpha
$
is $\succcurlyeq x_1$
and $\ne x_1$,
and yet it lies in $\tW(G,T)\,x_0$,
contradicting
the maximality assumption
of $x_1$.
Hence $x_1$ belongs to $C_0$.
Consequently,
there exists an element
$w\in \tW_0$
such that
$w(x_1)\in C(\Delta)$.
But $C(\Delta)$ is a fundamental domain
for the action of $\tW(G,T)$
on $\tX(T)_\bR$,
so $\tW(G,T)\,x_0$
meets $C(\Delta)$
at exactly one point,
namely $x_0$.
Since we also have
$w(x_1) \in C(\Delta) \cap \tW(G,T)\, x_0$,
it follows that
$w(x_1) = x_0$,
whence
$x_1 = w^{-1}(x_0) = x_0$
because $w^{-1} \in \tW_0$
fixes $x_0$. 
Thus $x_0 \succcurlyeq x$
as desired.
\end{proof}

\begin{lemma}
\label{wp edges:lemma3}
For every $\alpha\in\Delta(x_0)$,
there exists
an element $y\in\tY(T)_\bR$
such that
for every $\beta\in\Delta(x_0)$,
one has
$\langle \beta, y \rangle \geqslant 0$
and equality holds if and only if
$\beta = \alpha$.
\end{lemma}

\begin{proof}
First suppose that
the root $\alpha$
lies in $\Delta\bksl\Delta_0 \subseteq\Delta(x_0)$.
Choose an element
$\tilde{y}\in\tY(T)_\bR$
\textit{strongly dominant}
with respect to $\Delta$,
and set
$  y
   :=
   \tilde{y}
   -
   \tfrac{1}{2}\,
   \langle \alpha,\tilde{y} \rangle\,\alpha^\vee
   \in \tY(T)_\bR
$.
If $\beta\in\Delta(x_0)$,
then $\beta$ lies in $P(\Delta)$,
whence
\[
   \beta
   \ =\ 
   \sum_{\gamma\in\Delta}
      c_\gamma\,\gamma
   \qquad
   \text{for some coefficients
         $c_\gamma\in\bZ_{\geqslant 0}$,
	 $\gamma\in\Delta$}.
\]
It follows that
$  \langle \beta,y \rangle
   =
   \sum_{\gamma\in\Delta}
      c_\gamma\,A_\gamma,
$
where
$  A_\gamma
   :=
   \langle \gamma,y \rangle
   =
   \langle \gamma,\tilde{y} \rangle
   -
   \tfrac{1}{2}\,
   \langle \alpha,\tilde{y} \rangle
   \langle \gamma,\alpha^\vee\rangle
   \in \bR
$.
Now for every $\gamma\in\Delta$,
the two roots $\alpha,\gamma$
belong to the same
system $\Delta$ of simple roots,
so
$\langle \gamma,\alpha^\vee \rangle$
is $= 2$ if $\gamma=\alpha$
and is $\leqslant 0$ if $\gamma\ne\alpha$.
As $\langle \gamma,\tilde{y} \rangle$
is $>0$
for every $\gamma\in\Delta$
by assumption,
we infer that
$A_\gamma$ is $= 0$ if $\gamma = \alpha$
and is $> 0$ if $\gamma\ne\alpha$.
Consequently,
$\langle \beta,y \rangle$ is $\geqslant 0$
with equality if and only if
$\beta = c_\alpha\,\alpha$ for some $c_\alpha >0$,
and since the root datum of $(G,T)$ is reduced,
this can happen if and only if
$\beta = \alpha$.
This proves the lemma
in the case when
$\alpha \in \Delta\bksl\Delta_0
 \subseteq \Delta(x_0)$.

In general,
we can write
$\alpha = w_0(\alpha_0)$
for some
$w_0\in\tW_0$
and
$\alpha_0\in\Delta\bksl\Delta_0$.
By what we have already shown,
there exists an element $y_0\in\tY(T)_\bR$
such that
for every $\beta\in\Delta(x_0)$,
one has
$\langle \beta, y_0 \rangle \geqslant 0$
and equality holds if and only if
$\beta = \alpha_0$.
Set $y := w_0(y_0) \in \tY(T)_\bR$.
Then for any $\beta\in\Delta(x_0)$,
$ \langle \beta,y \rangle
 = \langle w_0^{-1}(\beta),y_0 \rangle $
is $\geqslant 0$ 
because $w_0^{-1}(\beta)$ belongs to $\Delta(x_0)$;
and equality holds
if and only if
$w_0^{-1}(\beta) = \alpha_0$,
which is to say
$\beta = \alpha$.
This proves the lemma.
\end{proof}

\subsection{}
Aside from
the results~(\ref{wp edges:lemma1}),
(\ref{wp edges:lemma2})
and
(\ref{wp edges:lemma3})
above,
the proof~(\ref{proof:thm:wp edges})
of theorem~(\ref{thm:wp edges})
also requires
the general
lemmas~(\ref{lemma:edge from vertex})
and~(\ref{lemma:face in half-ray})
below
about
the edges of a polytope
containing a given vertex.
We first recall
the basic facts
about polytopes
and their faces:

\smallskip
\noindent
\fact~A:
\ 
If $\polytope \subseteq V$ is
a polytope,
then
$ \polytope
  =
  \conv(\vt(\polytope))
$
is the convex hull
of its vertices.

\smallskip
\noindent
\fact~B:
\ 
If $F \subseteq \polytope$ is
a face of $\polytope$,
then
$F$ is also a polytope,
and
$ \vt(F)
  =
  F \cap \vt(\polytope)$.

\smallskip
\noindent
For the proof of
these standard results,
we refer to
\cite{Ziegler-LectPolytopes}
Prop.~2.2(i) and Prop.~2.3(i)
or \cite{Brondstead-IntroToConvexPolytopes}
Th.~7.2(c) and Th.~7.3.

\begin{lemma}
\label{lemma:edge from vertex}
Let $V$ be
a finite dimensional
$\bR$-vector space,
and let $\polytope \subseteq V$
be a polytope.
Let $x_0\in\vt(\polytope)$ be
a vertex.
A subset $e \subseteq \polytope$
of $\polytope$
is an edge containing $x_0$
if and only if:
\begin{itemize}
\item[(i)]
there exists a vertex
$x_1 \in \vt(\polytope)$
with $x_1 \ne x_0$,
such that $e = [x_0,x_1]$,
and
\item[(ii)]
there exist elements
$y \in V^\vee$
and
$c\in\bR$
such that
$ \langle x,y \rangle
  \leqslant
  c
$
for every vertex
$x \in \vt(\polytope)$,
and
equality holds
if and only if
$x = x_0$ or $x = x_1$.
\end{itemize}
\end{lemma}

\begin{proof}
First, suppose
$e\in\edg(\polytope)$ is
an edge containing $x_0$.
By \fact~B,
$e$ is a polytope
and $\vt(e) = e\cap\vt(\polytope)$.
Since $e$ has dimension~$1$,
$\vt(e)$ is necessarily
a set of two distinct elements,
namely
$x_0$ and another vertex
$x_1 \in \vt(\polytope)$
different from $x_0$.
By \fact~A,
we have
$e = [x_0,x_1]$,
which gives
condition~(i).
Condition~(ii)
follows from
the fact that
$e$ is
a face of $\polytope$
and that
$\{x_0,x_1\} = e \cap\vt(\polytope)$
from \fact~B.

Conversely,
suppose
$e \subseteq \polytope$
is a subset
satisfying~(i) and~(ii).
Since
$\polytope = \conv(\vt(\polytope))$
by \fact~A,
every element $x\in\polytope$
is a convex linear combination
of the points in $\vt(\polytope)$,
so
condition~(ii) implies that
$\langle x,y \rangle \leqslant c$
holds
for all $x \in \polytope$
and not just for
the vertices
$x \in \vt(\polytope)$.
Therefore,
\[
  F
  \ :=\ 
  \polytope
  \ \cap\ 
  \{\ 
     x\in V
     \ :\ 
     \langle x,y \rangle = c
  \ \}
  \qquad
  \text{is a face of $\polytope$}.
\]
Now
$\vt(F) = F \cap \vt(\polytope)$
by \fact~B,
and this intersection
is precisely $\{x_0,x_1\}$
according to condition~(ii).
So by \fact~A,
$ F = [x_0,x_1] = e $
is a face of $\polytope$
containing $x_0$;
since $x_0\ne x_1$
by condition~(i),
the face $e$
is of dimension~$1$.
\end{proof}

\begin{lemma}
\label{lemma:face in half-ray}
Let $V$ be
a finite dimensional
$\bR$-vector space,
and let $\polytope \subseteq V$
be a polytope.
Let $x_0\in\vt(\polytope)$ be
a vertex.
Suppose
there exist elements
$y \in V^\vee$
and
$c\in\bR$,
and an element
$v\in V$,
such that
$ \langle x,y \rangle
  \leqslant
  c
$
for every vertex
$x \in \vt(\polytope)$,
and equality holds
if and only if
$ x \in x_0 + \bR_{\geqslant 0}\,v$.
Let
$ F
  :=
  \polytope
  \ \cap\ 
  \{
      x \in V
      :
      \langle x,y \rangle = c
  \}
$.
Then:
\begin{itemize}
\item[(i)]
either $F$ is the singleton $\{x_0\}$,
\item[(ii)]
or there exists a vertex
$x_1\in\vt(\polytope)$ of $\polytope$
distinct from $x_0$
and lying in $F$,
in which case
$F$ is the line segment $[x_0,x_1]$
and is an edge of $\polytope$
containing $x_0$.
\end{itemize}
\end{lemma}

\begin{proof}
It follows from the hypotheses that
$F$ is a face of
the polytope $\polytope$.
From \fact~B,
we see that
$ \vt(F)
  =
  F
  \cap
  \vt(\polytope)
  =
  \{
      x \in \vt(\polytope)
      :
      \langle x,y \rangle = c
  \}
$,
which by assumption
is contained in
the half-line
$x_0 + \bR_{\geqslant 0}\,v$.
By \fact~A,
$F = \conv(\vt(F))$
is also contained in that half-line,
so it follows that
$F$ is of dimension~$0$ or $1$,
and correspondingly,
either
$\vt(F) = \{x_0\}$
or
$\vt(F) = \{x_0,x_1\}$
for some $x_1\in\vt(\polytope)$
with $x_1\ne x_0$;
these correspond to
cases~(i) and~(ii)
respectively.
\end{proof}

\subsection{Proof of~theorem~(\ref{thm:wp edges})}
\label{proof:thm:wp edges}
To show that
the map~(\ref{thm:wp edges:bijection})
is well-defined,
let $\alpha\in\Delta(x_0)$ be given;
we shall apply
lemma~(\ref{lemma:face in half-ray})
to show that
$[x_0,s_\alpha(x_0)]$
is an edge of $\polytope(\rho)$
containing $x_0$.
By~(\ref{wp edges:lemma3}),
there exists
an element $y\in\tY(T)_\bR$
such that
for every $\beta\in\Delta(x_0)$,
one has
$\langle \beta, y \rangle \geqslant 0$
and equality holds
if and only if
$\beta = \alpha$.
Set
$ c := \langle x_0,y \rangle \in\bR$.
If $w\in\tW(G,T)$ is
any element of the Weyl group,
then by~(\ref{wp edges:lemma2}),
\[
   w(x_0)
   \ =\ 
   x_0
   -
   \sum_{\beta\in\Delta(x_0)}
       c_\beta\,\beta
   \qquad
   \text{for some coefficients
         $c_\beta\in\bZ_{\geqslant 0}$,
	 $\beta\in\Delta(x_0)$}.
\]
It now follows from
the positivity property of
the element $y$
that
$\langle w(x_0),y \rangle$ is $\leqslant c$,
and equality holds
if and only if
$w(x_0) = x_0$ or
$w(x_0) = x_0 - c_\alpha\,\alpha$
for some $c_\alpha > 0$.
The hypotheses of
lemma~(\ref{lemma:face in half-ray})
are therefore satisfied,
and we infer that
\[
   F
   \ :=\ 
   \polytope(\rho)
   \ \cap\ 
   \{\ 
       x \in \tX(T)_\bR
       \ :\ 
       \langle x, y \rangle = c
   \ \}
\]
is either $\{x_0\}$
or an edge of $\polytope(\rho)$
containing $x_0$.
Since $\alpha\in\Delta(x_0)$
does not lie in $P_0$
by~(\ref{wp edges:lemma1}),
one has
$\langle x_0,\alpha^\vee \rangle > 0$,
so
$ s_\alpha(x_0)
  =
  x_0
  -
  \langle x_0,\alpha^\vee\rangle\,\alpha
$
is $\ne x_0$ and lies in $F$;
this rules out
the case $F = \{x_0\}$,
and thus
we must have
$F = [x_0,s_\alpha(x_0)]$,
and it belongs to
$\edges(\rho,{x_0})$.

To show that
the map~(\ref{thm:wp edges:bijection})
is injective,
suppose we are given
$\alpha,\beta\in \Delta(x_0)$
such that
$ [x_0,s_\alpha(x_0)]
  =
  [x_0,s_{\beta}(x_0)]$
in $\edges(\rho,{x_0})$;
then
$ s_\alpha(x_0)
  =
  s_{\beta}(x_0)
  \ne
  x_0 $,
from which we obtain
$ \langle x_0, \alpha^\vee \rangle \, \alpha
  =
  \langle x_0, {\beta}^\vee \rangle \, \beta
  \ne 0 $.
Since the root datum of $(G,T)$ is reduced,
and since $\alpha, \beta$ are positive roots
by~(\ref{wp edges:lemma1}),
it follows that
$\alpha = \beta$.

To show that
the map~(\ref{thm:wp edges:bijection})
is surjective,
let $e\in\edges(\rho,{x_0})$
be an edge containing $x_0$.
Applying
lemma~(\ref{lemma:edge from vertex}),
we see that
$e = [x_0, x_1]$
for some $x_1 \in\tW(G,T)\,x_0$
with $x_1\ne x_0$,
and there are elements
$y \in \tY(T)_\bR$ and $c\in\bR$
such that
for every $w\in\tW(G,T)$,
one has
$\langle w(x_0), y \rangle \leqslant c$,
and equality holds
if and only if
$w(x_0) = x_0$ or $w(x_0) = x_1$.
We claim that
$x_1 = s_\alpha(x_0)$
for some $\alpha\in \Delta(x_0)$.
Suppose this is not the case.
Then for each $\alpha\in \Delta(x_0)$,
we would have
$s_\alpha(x_0) \ne x_1$,
and since
$\alpha \notin P_0$
by~(\ref{wp edges:lemma1}),
we would also have
$s_\alpha(x_0) \ne x_0$,
and hence it follows that
$\langle s_\alpha(x_0),y \rangle < c$.
Since
$ s_\alpha(x_0)
  =
  x_0
  -
  \langle x_0,\alpha^\vee \rangle\,\alpha
$,
we infer that
$ \langle x_0,\alpha^\vee \rangle
  \,
  \langle \alpha,y \rangle$
is $> 0$,
and since $x_0$ is dominant
with respect to $\Delta$,
this implies
$\langle \alpha,y \rangle$
is $> 0$,
for each $\alpha\in \Delta(x_0)$.
However,
by~(\ref{wp edges:lemma2}),
we can write
\[
   x_1
   \ =\ 
   x_0
   -
   \sum_{\alpha\in \Delta(x_0)}
       c_\alpha\,\alpha
   \qquad
   \text{for some coefficients
         $c_\alpha \in \bZ_{\geqslant 0}$,
	 $\alpha\in \Delta(x_0)$},
\]
with at least one of
the $c_\alpha$'s being $>0$
because $x_1 \ne x_0$.
If we apply the pairing
$\langle -,y \rangle$
to both sides of the equation,
we see that
the sum
$ \sum_{\alpha\in \Delta(x_0)}
      c_\alpha\,\langle \alpha,y \rangle
$
is~$0$,
which is a contradiction
since this is a sum of non-negative terms
with at least one strictly positive term.
This proves our claim,
and hence
$ e = [x_0,s_\alpha(x_0)] $
for some $\alpha\in\Delta(x_0)$.
\qed


\subsection*{Acknowledgments}
I am grateful to
Wee Teck Gan,
Johan de~Jong,
and Nick Katz
for useful discussions
during the course of this work.
%




\end{document}